\documentclass{elsart3-1}

\usepackage[utf8x]{inputenc}
\usepackage[english,francais]{babel}
\usepackage{amsmath,amsfonts,amssymb}
\usepackage{subfigure}
\usepackage{cite}
\usepackage{multicol}
\usepackage{tikz}

\newtheorem{Theoreme}{Th\'eor\`eme}[section]
\newtheorem{Proposition}[Theoreme]{Proposition}
\renewcommand{\theequation}{\arabic{equation}}
\renewcommand{\leq}{\leqslant}
\renewcommand{\geq}{\geqslant}

\newcommand{\K}{\mathbb{K}}
\newcommand{\EnsPermu}{\mathfrak{S}}
\newcommand{\EnsNat}{\mathbb{N}}
\newcommand{\EnsRel}{\mathbb{Z}}
\newcommand{\CalP}{\mathcal{P}}
\newcommand{\CalB}{\mathcal{B}}
\newcommand{\T}{{\sf T}}
\newcommand{\TT}{{\it T}}

\newcommand{\End}{{\sf End}}
\newcommand{\FP}{{\sf FP}}
\newcommand{\MT}{{\sf MT}}
\newcommand{\APE}{{\sf APE}}
\newcommand{\FCat}[1]{{\sf FCat^{(#1)}}}
\newcommand{\Schr}{{\sf Schr}}
\newcommand{\Motz}{{\sf Motz}}
\newcommand{\Comp}{{\sf Comp}}
\newcommand{\AnD}{{\sf AnD}}
\newcommand{\SComp}{{\sf SComp}}
\newcommand{\DD}{{\sf D}}
\newcommand{\La}{{\tt \textcolor{Rouge}{a}}}
\newcommand{\Lb}{{\tt \textcolor{Bleu}{b}}}
\newcommand{\Lc}{{\tt \textcolor{Vert}{c}}}

\definecolor{Noir}{RGB}{0,0,0}
\definecolor{Rouge}{RGB}{205,35,38}
\definecolor{Bleu}{RGB}{2,60,195}
\definecolor{Vert}{RGB}{63,143,1}
\definecolor{Violet}{RGB}{181,18,225}
\definecolor{Orange}{RGB}{255,113,15}

\tikzstyle{Injection} = [Noir!100,draw,>->]
\tikzstyle{Surjection} = [Noir!100,draw,->>]

\begin{document}

\centerline{}
\begin{frontmatter}

\selectlanguage{francais}
\title{Construction d'op\'erades ensemblistes \`a partir de mono\"ides}
\journal{the Acad\'emie des sciences}
\author[SG]{Samuele Giraudo}
\ead{samuele.giraudo@univ-mlv.fr}
\address[SG]{Institut Gaspard Monge, Universit\'e Paris-Est
Marne-la-Vall\'ee, 5 Boulevard Descartes, Champs-sur-Marne,
77454 Marne-la-Vall\'ee cedex 2, France}

\medskip
\begin{center}
    {\small Re\c{c}u le *****~; accept\'e apr\`es r\'evision le ***** \\
    Pr\'esent\'e par *****}
\end{center}

\begin{abstract}
    \selectlanguage{francais}
    Nous étudions une construction fonctorielle de la catégorie des monoïdes
    vers la catégorie des opérades ensemblistes et donnons des exemples
    combinatoires d'applications. \\
    \noindent
    {\it Pour citer cet article~: S. Giraudo, C. R. Acad. Sci. Paris, Ser. * *** (****).}
    \vskip 0.5\baselineskip

    \selectlanguage{english}
    \noindent{\bf Abstract --- Constructing set-operads from monoids} \\
    \vskip 0.5\baselineskip
    \noindent
    We study a functorial construction from the category of monoids
    to the category of set-operads and we give some combinatorial
    examples of applications. \\
    \noindent
    {\it To cite this article: S. Giraudo, C. R. Acad. Sci. Paris, Ser. * *** (****).}
\end{abstract}

\end{frontmatter}

\selectlanguage{francais}

\vspace{-1.5em}
\section{Introduction\vspace{-1em}} \label{sec:Introduction}
Les \emph{opérades} sont des structures algébriques qui formalisent la
notion de composition d'opérateurs et les relations qu'ils vérifient. Plus
précisément, une opérade contient des opérateurs munis de~$n \geq 1$ entrées
et d'une unique sortie. Deux opérateurs~$x$ et~$y$ peuvent être composés
en~$i\ieme$ position en greffant la sortie de~$y$ sur la~$i\ieme$ entrée
de~$x$. Le nouvel opérateur ainsi obtenu est noté~$x \circ_i y$. Il est
de plus possible dans une opérade de permuter les entrées d'un opérateur~$x$
en faisant agir une permutation~$\sigma$. Le nouvel opérateur ainsi obtenu
est noté~$x \cdot \sigma$. L'un des principaux points forts de cette théorie
est qu'elle offre un cadre et un formalisme général pour étudier de
manière unifiée différents types d'algèbres, comme les algèbres associatives
et les algèbres de Lie. Dans cette article, nous considérons exclusivement
les \emph{opérades ensemblistes} qui sont des ensembles de la
forme~$\CalP := \biguplus_{n \geq 1} \CalP(n)$ où les~$\CalP(n)$ sont des
ensembles d'éléments d'\emph{arité}~$n$, munis d'\emph{applications de greffe}
\begin{equation}
    \circ_i : \CalP(n) \times \CalP(m) \to \CalP(n + m - 1),
    \qquad n, m \geq 1 \mbox{ et } 1 \leq i \leq n,
\end{equation}
et d'une \emph{action du groupe symétrique}
\begin{equation}
    \cdot : \CalP(n) \times \EnsPermu_n \to \CalP(n), \qquad n \geq 1,
\end{equation}
qui vérifient des axiomes naturels.

Nous proposons dans ce travail une construction fonctorielle~$\T$ qui
permet d'obtenir des opérades~$\T M$ à partir de monoïdes~$M$. Les éléments
de~$\T M$ d'arité~$n$ sont les mots de longueur~$n$ sur~$M$ vu comme un
alphabet, et l'expression de la greffe dans cette opérade s'obtient
directement par l'expression du produit de~$M$.

Dans des travaux antérieurs, Berger et Moerdijk~\cite{BM03} proposèrent
une construction~$\TT$ qui permet d'obtenir, à partir d'une bigèbre
commutative~$\CalB$, une coopérade~$\TT \CalB$. Notre construction~$\T$
et la construction~$\TT$ de ces deux auteurs sont différentes mais coïncident
dans de nombreux cas. Par exemple, lorsque~$(M, \bullet)$ est un monoïde
tel que pour tout~$x \in M$, l'ensemble des couples~$(y, z) \in M^2$ qui
vérifient~$y \bullet z = x$ est fini, alors l'opérade~$\T M$ est la duale
de la coopérade~$\TT \CalB$ où~$\CalB$ est la bigèbre duale de la bigèbre~$\K [M]$
munie du coproduit diagonal ($\K$ est un corps). En revanche, il existe des
opérades que l'on peut obtenir par la construction~$\T$ mais pas par la
construction~$\TT$ --- et réciproquement. Par exemple, l'opérade~$\T \EnsRel$
où~$\EnsRel$ est le monoïde additif des entiers relatifs ne peut être
obtenue comme duale d'une coopérade constructible par la construction de
Berger et Moerdijk.

En outre, notre construction est définie dans la catégorie des ensembles
et les calculs y sont explicites. Il est donc possible, à partir d'un
monoïde~$M$ quelconque de calculer simplement, si nécessaire à l'aide de
l'ordinateur, dans l'opérade~$\T M$.

Dans cet article, nous étudions plusieurs applications de la construction~$\T$
et mettons l'accent sur son caractère combinatoire. Plus précisément, nous
définissons, à partir de monoïdes usuels --- comme le monoïde additif des
entiers naturels ou les monoïdes cycliques --- diverses opérades qui mettent
en jeu plusieurs objets combinatoires connus. Nous construisons ainsi des
opérades sur des objets qui n'étaient pas pourvus d'une telle structure~:
arbres d'arité fixée, chemins de Motzkin, compositions d'entiers, animaux
dirigés et compositions d'entiers segmentées. Nous obtenons aussi de nouvelles
opérades sur des objets déjà pourvus d'une telle structure~: fonctions de
parking, mots tassés, arbres plans enracinés et arbres de Schröder. Notre
construction permet également de retrouver des opérades déjà connues par
ailleurs, comme l'opérade magmatique, l'opérade commutative associative
et l'opérade diassociative~\cite{Lod01}.

\vspace{-1.5em}
\section{Un foncteur des monoïdes vers les opérades ensemblistes\vspace{-1em}} \label{sec:Foncteur}
Soit~$(M, \bullet)$ un monoïde. Définissons~$\T M$ comme
l'ensemble~$\T M := \biguplus_{n \geq 1} \T M(n)$, où pour tout~$n \geq 1$,
\begin{equation}
    \T M(n) := \left\{(x_1, \dots, x_n) :
        x_i \in M \mbox{ pour tout $1 \leq i \leq n$}\right\}.
\end{equation}
Les éléments de~$\T M(n)$ sont ainsi les mots sur l'alphabet~$M$ de
longueur~$n$. Munissons maintenant l'ensemble~$\T M$ d'applications de greffe
\begin{equation} \label{eq:TDomaineSubs}
    \circ_i : \T M(n) \times \T M(m) \to \T M(n + m - 1),
    \qquad n, m \geq 1 \mbox{ et } 1 \leq i \leq n,
\end{equation}
définies pour tous~$x \in \T M(n)$, $y \in \T M(m)$ et~$1 \leq i \leq n$ par
\begin{equation} \label{eq:TSub}
    x \circ_i y :=
    (x_1, \; \dots, \; x_{i-1},
    \; x_i \bullet y_1, \; \dots, \; x_i \bullet y_m,
    \; x_{i+1}, \; \dots, \; x_n).
\end{equation}
Par exemple, si $M$ est le monoïde additif des entiers
naturels, nous avons dans~$\T M$,
$\textcolor{Bleu}{2}{\bf 1} \textcolor{Bleu}{23} \circ_2 \textcolor{Rouge}{30313}
= \textcolor{Bleu}{2}\textcolor{Rouge}{41424}\textcolor{Bleu}{23}$.
Munissons également chaque ensemble~$\T M(n)$ d'une action à droite du
groupe symétrique
\begin{equation}
    \cdot : \T M(n) \times \EnsPermu_n \to \T M(n), \qquad n \geq 1,
\end{equation}
définie pour tous~$x \in \T M(n)$ et~$\sigma \in \EnsPermu_n$ par
\begin{equation}
    x \cdot \sigma := \left(x_{\sigma_1}, \dots, x_{\sigma_n}\right).
\end{equation}
Par exemple, si~$\Lb\Lb\Lc\Lb\La$ est un élément de $\T M$ et~$\sigma$ est
la permutation~$23514$, nous avons~$\Lb\Lb\Lc\Lb\La \cdot \sigma = \Lb\Lc\La\Lb\Lb$.
\smallskip

Si~$M$ et~$N$ sont deux monoïdes et~$\theta : M \to N$ un morphisme de
monoïdes, notons~$\T \theta$ l'application
\begin{equation}
    \T \theta : \T M \to \T N,
\end{equation}
définie pour tout~$(x_1, \dots, x_n) \in \T M(n)$ par
\begin{equation}
    \T \theta\left(x_1, \dots, x_n\right) :=
    \left(\theta(x_1), \dots, \theta(x_n)\right).
\end{equation}

\begin{Theoreme}
    La construction~$\T$ est un foncteur de la catégorie des monoïdes avec
    morphismes de monoïdes vers la catégorie des opérades ensemblistes
    avec morphismes d'opérades ensemblistes. De plus,~$\T$ respecte les
    injections et les surjections.
\end{Theoreme}

\vspace{-1.5em}
\section{Quelques opérades obtenues par le foncteur~$\T$\vspace{-1em}} \label{sec:Exemples}
Pour illustrer la richesse combinatoire de la construction~$\T$, nous
construisons à présent des sous-opérades symétriques ou non de l'opérade
obtenue à partir du monoïde additif des entiers naturels~$\EnsNat$ et
énonçons quelques-unes de leurs propriétés. Dans ce qui suit,~$\EnsNat_2$
(resp.~$\EnsNat_3$) désigne le monoïde des entiers naturels modulo~$2$
(resp. $3$). La figure~1 répertorie les relations entre ces opérades.
\begin{figure}[ht]
    \centering
    \parbox{.5\textwidth}{
    \scalebox{.65}{
    \begin{tikzpicture}[scale=.6]
        \node(TN)at(2,0){$\T \EnsNat$};
        \node(TN2)at(-4,-2){$\T \EnsNat_2$};
        \node(TN3)at(8,-2){$\T \EnsNat_3$};
        \node(End)at(-.5,-2){$\End$};
        \node(FP)at(-1.25,-4){$\FP$};
        \node(MT)at(-1.25,-6){$\MT$};
        \node(Schr)at(-.5,-8){$\Schr$};
        \node(FCat1)at(2,-10){$\FCat{1}$};
        \node(FCat2)at(2,-8){$\FCat{2}$};
        \node(FCat3)at(2,-6){$\FCat{3}$};
        \node(SComp)at(8,-10){$\SComp$};
        \node(AnD)at(8,-12){$\AnD$};
        \node(APE)at(5,-12){$\APE$};
        \node(Motz)at(-1,-12){$\Motz$};
        \node(Comp)at(-4,-12){$\Comp$};
        \node(FCat0)at(2,-14){$\FCat{0}$};
        \draw[Surjection](TN)--(TN2);
        \draw[Surjection](TN)--(TN3);
        \draw[Injection](End)--(TN);
        \draw[Injection](FP)--(End);
        \draw[Injection](MT)--(FP);
        \draw[Injection](Schr)--(MT);
        \draw[Injection](FCat1)--(Schr);
        \draw[Injection](FCat1)--(FCat2);
        \draw[Injection](FCat2)--(FCat3);
        \draw[Injection,dashed](FCat3)--(TN);
        \draw[Surjection](FCat2)--(SComp);
        \draw[Injection](SComp)--(TN3);
        \draw[Surjection](FCat1)--(AnD);
        \draw[Injection](AnD)--(SComp);
        \draw[Injection](APE)--(FCat1);
        \draw[Injection](Motz)--(FCat1);
        \draw[Surjection](FCat1)--(Comp);
        \draw[Injection](Comp)--(TN2);
        \draw[Surjection](Comp)--(FCat0);
        \draw[Surjection](APE)--(FCat0);
        \draw[Surjection](AnD)--(FCat0);
        \draw[Injection](FCat0)--(Motz);
    \end{tikzpicture}}}
    \parbox{.4\textwidth}{
    \caption{Le diagramme des sous-opérades et quotients non-symétriques de
    l'opérade~$\T \EnsNat$. Les flèches~$\rightarrowtail$
    (resp.~$\twoheadrightarrow$) sont des morphismes injectifs (resp.
    surjectifs) d'opérades non-symétriques.}}
\end{figure}
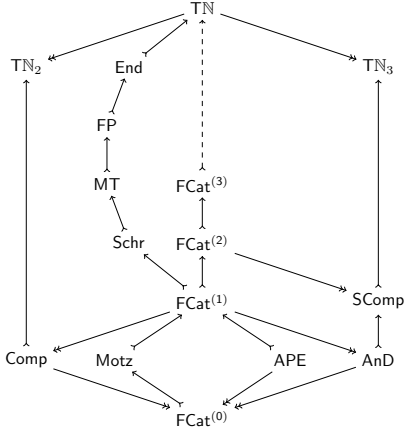

Un mot~$u$ est une endofonction (resp. fonction de parking, mot tassé)
\emph{tordue} si le mot $(u_1 + 1, u_2 + 1, \dots, u_{|u|} + 1)$
est une endofonction (resp. fonction de parking, mot tassé). Notons~$\End$
(resp.~$\FP$, $\MT$) l'ensemble des endofonctions (resp. fonctions de parking,
mots tassés).
\smallskip

\begin{Proposition}
    Les ensembles~$\End$, $\FP$ et~$\MT$ forment des sous-opérades
    symétriques de~$\T \EnsNat$. De plus,~$\MT$ est engendrée, en tant
    qu'opérade symétrique, par~$00$ et~$01$.
\end{Proposition}
\smallskip

\begin{Theoreme}
    Soit~$\APE$ la sous-opérade non-symétrique de~$\T \EnsNat$ engendrée
    par~$01$. Alors, les éléments de~$\APE$ d'arité~$n$ sont exactement
    les mots~$x$ sur l'alphabet~$\EnsNat$ qui vérifient~$x_1 = 0$
    et~$1 \leq x_{i + 1} \leq x_i + 1$ pour tout~$1 \leq i \leq n - 1$.
    De plus, les éléments de~$\APE$ d'arité~$n$ sont en bijection avec
    les arbres plans enracinés à~$n$ n\oe uds. Enfin,~$\APE$ est isomorphe
    à l'opérade non-symétrique libre sur un générateur d'arité deux.
\end{Theoreme}
\smallskip

\begin{Theoreme}
    Soit~$k \geq 0$ et~$\FCat{k}$ la sous-opérade non-symétrique
    de~$\T \EnsNat$ engendrée par~$00$, $01$, \dots, $0k$. Alors, les
    éléments de~$\FCat{k}$ d'arité~$n$ sont exactement les mots~$x$ sur
    l'alphabet~$\EnsNat$ qui vérifient $x_1 = 0$ et~$0 \leq x_{i + 1} \leq x_i + k$
    pour tout~$1 \leq i \leq n - 1$. De plus, les éléments de~$\FCat{k}$
    d'arité~$n$ sont en bijection avec les arbres plans enracinés
    d'arité~$k + 1$ et de taille~$n$. Enfin,~$\FCat{2}$ est isomorphe à
    l'opérade non-symétrique libre engendrée par deux générateurs~$\La$
    et~$\Lb$ d'arité deux, sujets aux trois relations
    \vspace{-2.7em}
    \begin{multicols}{3}
        \begin{equation}
            \La \circ_1 \La = \La \circ_2 \La,
        \end{equation}

        \begin{equation}
            \Lb \circ_1 \La = \La \circ_2 \Lb,
        \end{equation}

        \begin{equation}
            \Lb \circ_1 \Lb = \Lb \circ_2 \La.
        \end{equation}
    \end{multicols}
\end{Theoreme}

\begin{Theoreme}
    Soit~$\Schr$ la sous-opérade non-symétrique de~$\T \EnsNat$ engendrée
    par~$00$, $01$ et~$10$. Alors, les éléments de~$\Schr$ sont exactement
    les mots~$x$ sur l'alphabet~$\EnsNat$ qui ont au moins une occurrence
    de~$0$ et tels que pour toute lettre~$b \geq 1$ de~$x$, il existe
    une lettre~$a = b - 1$ telle que~$x$ possède un facteur~$a u b$
    ou~$b u a$ où~$u$ est un mot composé de lettres~$c$ vérifiant~$c \geq b$.
    De plus, les éléments de~$\Schr$ d'arité~$n$ sont en bijection avec
    les arbres de Schröder~\cite{FS09} à~$n$ feuilles. Enfin,~$\Schr$
    est isomorphe à l'opérade non-symétrique libre engendrée par trois
    générateurs~$\La$, $\Lb$ et~$\Lc$ d'arité deux, sujets aux sept relations
    \vspace{-2.7em}
    \begin{multicols}{4}
        \begin{equation}
            \Lb \circ_1 \La = \La \circ_2 \Lb,
        \end{equation}
        \begin{equation}
            \Lb \circ_1 \Lb = \Lb \circ_2 \La,
        \end{equation}

        \begin{equation}
            \La \circ_1 \Lc = \Lc \circ_2 \La,
        \end{equation}
        \begin{equation}
            \Lc \circ_1 \La = \Lc \circ_2 \Lc,
        \end{equation}

        \begin{equation}
            \La \circ_1 \La = \La \circ_2 \La,
        \end{equation}
        \begin{equation}
            \Lb \circ_1 \Lc = \Lc \circ_2 \Lb,
        \end{equation}

        \begin{equation}
            \La \circ_1 \Lb = \La \circ_2 \Lc.
        \end{equation}
    \end{multicols}
\end{Theoreme}

\begin{Theoreme}
    Soit~$\Motz$ la sous-opérade non-symétrique de~$\T \EnsNat$ engendrée
    par~$00$ et~$010$. Alors, les éléments de~$\Motz$ d'arité~$n$ sont
    exactement les mots sur l'alphabet~$\EnsNat$ qui commencent et se terminent
    par~$0$ et tels que~$|x_i - x_{i + 1}| \leq 1$ pour tout~$1 \leq i \leq n - 1$.
    De plus, les éléments de~$\Motz$ d'arité~$n$ sont en bijection avec
    les chemins de Motzkin~\cite{FS09} à~$n$ pas. Enfin,~$\Motz$ est
    isomorphe à l'opérade non-symétrique libre engendrée par un
    générateur~$\La$ d'arité deux et un générateur~$\Lb$ d'arité trois,
    sujets aux quatre relations
    \vspace{-2.7em}
    \begin{multicols}{4}
        \begin{equation}
            \La \circ_1 \La = \La \circ_2 \La,
        \end{equation}

        \begin{equation}
            \La \circ_1 \Lb = \Lb \circ_3 \La,
        \end{equation}

        \begin{equation}
            \Lb \circ_1 \La = \La \circ_2 \Lb,
        \end{equation}

        \begin{equation}
            \Lb \circ_1 \Lb = \Lb \circ_3 \Lb.
        \end{equation}
    \end{multicols}
\end{Theoreme}

\begin{Theoreme}
    Soit~$\Comp$ la sous-opérade non-symétrique de~$\T \EnsNat_2$ engendrée
    par~$00$ et~$01$. Alors, les éléments de~$\Comp$ sont exactement les
    mots sur l'alphabet~$\{0, 1\}$ qui commencent par~$0$. De plus, les
    éléments de~$\Comp$ d'arité~$n$ sont en bijection avec les compositions
    de l'entier~$n$. Enfin,~$\Comp$ est isomorphe à l'opérade non-symétrique
    libre engendrée par deux générateurs~$\La$ et~$\Lb$ d'arité deux,
    sujets aux quatre relations
    \vspace{-2.7em}
    \begin{multicols}{4}
        \begin{equation}
            \La \circ_1 \La = \La \circ_2 \La,
        \end{equation}

        \begin{equation}
            \Lb \circ_1 \La = \La \circ_2 \Lb,
        \end{equation}

        \begin{equation}
            \Lb \circ_1 \Lb = \Lb \circ_2 \La,
        \end{equation}

        \begin{equation}
            \La \circ_1 \Lb = \Lb \circ_2 \Lb.
        \end{equation}
    \end{multicols}
\end{Theoreme}

\begin{Proposition}
    Soit~$\AnD$ la sous-opérade non-symétrique de~$\T \EnsNat_3$ engendrée
    par~$00$ et~$01$. Alors, les éléments de~$\AnD$ d'arité~$n$ sont
    en bijection avec les animaux dirigés~\cite{FS09} de taille~$n$.
    De plus,~$\AnD$ n'admet pas de présentation quadratique.
\end{Proposition}
\smallskip

\begin{Theoreme}
    Soit~$\SComp$ la sous-opérade non-symétrique de~$\T \EnsNat_3$ engendrée
    par~$00$, $01$ et~$02$. Alors, les éléments de~$\SComp$ sont exactement
    les mots sur l'alphabet~$\{0, 1, 2\}$ qui commencent par~$0$. De plus,
    les éléments de~$\SComp$ d'arité~$n$ sont en bijection avec les
    compositions segmentées~\cite{FS09} de l'entier~$n$. Enfin,~$\SComp$
    est isomorphe à l'opérade non-symétrique libre engendrée par trois
    générateurs~$\La$, $\Lb$ et~$\Lc$ d'arité deux, sujets aux neuf relations
    \vspace{-2.7em}
    \begin{multicols}{3}
        \begin{equation}
            \La \circ_1 \La = \La \circ_2 \La,
        \end{equation}
        \begin{equation}
            \Lb \circ_1 \La = \La \circ_2 \Lb,
        \end{equation}
        \begin{equation}
            \Lb \circ_1 \Lb = \Lb \circ_2 \La,
        \end{equation}

        \begin{equation}
            \Lc \circ_1 \La = \La \circ_2 \Lc,
        \end{equation}
        \begin{equation}
            \Lc \circ_1 \Lc = \Lc \circ_2 \La,
        \end{equation}
        \begin{equation}
            \Lb \circ_1 \Lc = \Lc \circ_2 \Lc,
        \end{equation}

        \begin{equation}
            \Lc \circ_1 \Lb = \Lb \circ_2 \Lb,
        \end{equation}
        \begin{equation}
            \La \circ_1 \Lb = \Lb \circ_2 \Lc,
        \end{equation}
        \begin{equation}
            \La \circ_1 \Lc = \Lc \circ_2 \Lb.
        \end{equation}
    \end{multicols}
\end{Theoreme}

\begin{Proposition}
    Soit le monoïde~$M := \{0, 1\}$ muni de la multiplication des entiers
    comme produit. Soit~$\DD$ la sous-opérade de~$\T M$ engendrée par~$01$
    et~$10$. Alors, les éléments de~$\DD$ sont exactement les mots qui
    contiennent exactement une occurrence de~$1$. De plus,~$\DD$ est isomorphe
    à l'opérade diassociative~\cite{Lod01}.
\end{Proposition}

\vspace{-1.5em}
\renewcommand{\refname}{Références\vspace{-1em}}
\bibliographystyle{plain}
\bibliography{Bibliographie}

\end{document}